\documentclass[12pt]{article}

\usepackage{amsmath,amssymb,url}
\usepackage{tikz,color}

\newtheorem{theorem}{Theorem}[section]
\newtheorem{corollary}[theorem]{Corollary}

\newtheorem{proposition}[theorem]{Proposition}

\newtheorem{problem}[theorem]{Problem}

\newtheorem{remark}[theorem]{Remark}

\newcommand{\ecc}{{\rm ec}}
\newcommand{\rad}{{\rm rad}}
\newcommand{\diam}{{\rm diam}}
\newcommand{\tr}{{\rm Tr}}
\newcommand{\s}{\color{black}}

\def\cp{\,\square\,}
\newcommand{\proof}{\noindent{\bf Proof.\ }}
\newcommand{\qed}{\hfill $\square$ \bigskip}

\vspace{8mm}
\title{Comparing Wiener complexity with eccentric complexity}

\author{Kexiang Xu$\/^{a}$, Aleksandar Ili\'{c}$\/^{b}$, Vesna Ir\v{s}i\v{c}$\/^{c,d}$, Sandi Klav\v{z}ar$^{c,d,e}$, \\
Huimin Li$\/^{a}$ \\\\
 $^{a}$ \small  College of Science, Nanjing University of
 Aeronautics \& Astronautics,\\
 \small Nanjing, Jiangsu 210016, PR China\\
$^{b}$ \small Facebook Inc, Menlo Park, California 94025, USA\\
$^{c}$ \small Institute of Mathematics, Physics and Mechanics, Ljubljana, Slovenia \\
 $^{d}$ \small Faculty of Mathematics and Physics, University of Ljubljana, Slovenia\\
  $^{e}$ \small Faculty of Natural Sciences and Mathematics, University of Maribor, Slovenia\\
\small {\tt kexxu1221@126.com} (K.\ Xu)\\
\small {\tt aleksandari@gmail.com} (A.\ Ili\'{c})\\
\small {\tt vesna.irsic@fmf.uni-lj.si} (V.\ Ir\v{s}i\v{c})\\
\small {\tt sandi.klavzar@fmf.uni-lj.si} (S.\ Klav\v{z}ar)\\
\small {\tt 1213966965@qq.com} (H.\ Li)
}
\date{}

\begin{document}

\maketitle


\begin{abstract}
The transmission of a vertex $v$ of a graph $G$ is the sum of distances from $v$ to all the other vertices in $G$. The Wiener complexity of $G$ is the number of different transmissions of its vertices. Similarly, the eccentric complexity of $G$ is defined as the number of different eccentricities of its vertices. In this paper these two complexities are compared. The complexities are first studied on Cartesian product graphs.  Transmission indivisible graphs and arithmetic transmission graphs are introduced to demonstrate sharpness of upper and lower bounds on the Wiener complexity, respectively. It is shown that for almost all graphs the Wiener complexity is not smaller than the eccentric complexity.  This property is proved for trees, the equality holding precisely for center-regular trees. Several families of graphs in which the complexities are equal are constructed. Using the Cartesian product, it is proved that the eccentric complexity can be arbitrarily larger than the Wiener complexity.  Additional infinite families of graphs with this property are constructed by amalgamating  universally diametrical graphs with center-regular trees.
\end{abstract}


\medskip\noindent
\textbf{AMS Math.\ Subj.\ Class.\ (2010)}: 05C12, 05C76

\section{Introduction}

If $G = (V(G), E(G))$ is a graph, then $d_G(u,v)$ denotes the shortest-path distance between vertices $u, v\in V(G)$. The {\em transmission} $\tr_G(v)$ of a vertex $v\in V(G)$ is the sum of distances from $v$ to the vertices in $G$, that is,
$$\tr_G(v)=\sum\limits_{u\in V(G)}d_G(u,v)\,.$$ The Wiener index $W(G)$ of $G$ can then be defined as $W(G)=\frac{1}{2}\sum\limits_{v\in V(G)}\tr_G(v)$. The {\em eccentricity} $\ecc_G(v)$ of a vertex $v\in V(G)$ is the maximum distance from $v$ to other vertices of $G$:
$$\ecc_G(v)=\max\limits_{u\neq v}d_G(u,v)\,.$$
Eccentricity is a central concept of metric graph theory and has many applications elsewhere, in particular in location theory and in chemical graph theory. In the latter area, important eccentricity-based graph invariants (alias topological indices in mathematical chemistry) include the first and the second Zagreb eccentricity indices~\cite{VG2010}, eccentric connectivity index~\cite{SGM1997, zhang-2019}, and connective eccentricity index~\cite{GSM2000}. For mathematical properties of these invariants see~\cite{DZT2012, QZL2017, XAD2017, XDL2016, XuLi2016}.

The \textit{Wiener complexity} $C_W(G)$ of a graph $G$ was introduced in~\cite{AAKS2014} (under the name Wiener dimension) as the number of different transmission of vertices in $G$:
$$C_W(G) = |\{ \tr_G(v):\ v\in V(G) \}|\,.$$
The Wiener complexity of graphs has been further investigated in~\cite{AK2018, jamilet-2016, KJRMP2018+}. In the same spirit as the Wiener complexity is defined, the connective eccentric complexity~\cite{AK2016} and the eccentric complexity~\cite{ADX2017+} have been recently introduced. The {\em eccentric complexity} $C_{\ecc}(G)$ of a graph $G$ is the number of different eccentricities in $G$. Equivalently,
\begin{equation}
\label{eq:C_ec}
C_{\ecc}(G) = \diam(G) - \rad(G) + 1\,,
\end{equation}
where $\diam(G) = \max\limits_{v\in V(G)}\ecc_G(v)$ is the \textit{diameter} of $G$ and $\rad(G) = \min\limits_{v\in V(G)}\ecc_G(v)$ is the \textit{radius} of $G$.

In view of the conceptional similarities between the Wiener complexity and the eccentric complexity, we compare in this paper these two complexities and proceed as follows. In the rest of this section we list definitions, concepts, and known results needed. In Section~\ref{sec:Cartesian} we consider the Wiener complexity and the eccentric complexity of Cartesian products.  In particular, a new lower bound on the Wiener complexity is proved and shown to be sharp using the so-called arithmetic transmission graphs. To demonstrate that the Wiener complexity of a Cartesian product can be equal to the product of the complexities of the factors, transmission indivisible graphs are introduced. In Section~\ref{sec:compare} we first prove that $C_{\ecc}(G) \le C_W(G)$ holds for almost all graphs $G$. Consequently, in Subsections~\ref{subsec:equality-case} and~\ref{subsec:larger}, we consider the graphs $G$ for which $C_{\ecc}(G) = C_W(G)$ and $C_{\ecc}(G) > C_W(G)$ holds, respectively.   We prove that for a tree $T$ we always have $C_{\ecc}(G) \le C_W(G)$, and that equality holds precisely for center-regular trees. We construct several families of graphs, among them two families of product graphs, for which the equality holds. Using the Cartesian product, we find an infinite family of graphs $G$ with the property $C_{\ecc}(G) > C_W(G)$. Moreover, the construction shows that the difference can be arbitrarily large. Finally,  amalgamating universally diametrical graphs with center-regular trees we construct additional infinite families of graphs $G$ for which  $C_{\ecc}(G) > C_W(G)$ holds.

\subsection{Preliminaries}

If $k$ is a positive integer, then $[k]=\{1, \ldots, k\}$. The degree of a vertex $v$ of a graph $G$ is denoted by $\deg_G(v)$. If $G$ is a graph, then $n(G)$ denotes the order of $G$.

The \emph{center} $C(G)$ of a graph $G$ is the set of vertices of $G$ of the minimum eccentricity, these vertices being called {\em central}. $G$ is \textit{self-centered}~\cite{buckley-1989} if all its vertices have the same eccentricity, that is, if and only if $C_{\ecc}(G)=1$.

Let $G$ be a graph. The \textit{transmission set} $\tr(G)$ of $G$ is the set of the transmissions of its vertices, that is, $\tr(G) = \{\tr_G(v):\ v\in V(G)\}$. The \textit{eccentricity set} ${\rm Ec}(G)$ of a graph $G$ is the set of the eccentricities of its vertices, that is, ${\rm Ec}(G) = \{\ecc_G(v):\ v\in V(G)\}$. A graph $G$ is \textit{transmission regular}~\cite{LDW2016} if all its vertices have the same transmission. In other words, transmission regular graphs are precisely the graphs $G$ with $C_W(G)=1$. In addition, $G$ is \textit{transmission irregular}~\cite{AK2018} if all its vertices have pairwise different transmissions, that is, if and only if $C_W(G)=n(G)$. We will make use of the following easy result on the transmission.

\begin{proposition}{\rm (\cite{AAKS2014})} \label{prp:2-ec}
Let $G$ be a graph with $\diam(G) = 2$. If $v\in V(G)$ with $\ecc_G(v)=2$, then $\tr_G(v)=2n(G)-2-\deg_G(v)$.
\end{proposition}

The {\em Cartesian product} $G\cp H$ of graphs $G$ and $H$ is the graph with $V(G\cp H) = V(G)\times V(H)$, vertices $(g,h)$ and $(g',h')$ being adjacent if $gg'\in E(G)$ and $h=h'$, or $g=g'$ and $hh'\in E(H)$.

\section{Complexities on Cartesian products}
\label{sec:Cartesian}

While comparing the Wiener complexity and the eccentric complexity we will extensively use the Cartesian product operation. In this section we hence recall known, and derive new related results.

It is well known that the distance function is additive on Cartesian product graphs. More precisely, if $G$ and $H$ are graphs, then
\begin{equation}
\label{eq:distance-formula}
d_{G\cp H}((g,h),(g^{\prime},h^{\prime}))=d_G(g,g^{\prime})+ d_H(h,h^{\prime})
\end{equation}
holds for arbitrary vertices $(g,h),\, (g^{\prime},h^{\prime})\in V(G\cp H)$, cf.~\cite [Proposition 5.1] {HIK2011}. This, in particular, implies that the diameter and radius are additive functions on Cartesian product graphs, which in turn gives the following closed formula for the eccentric complexity of Cartesian products.

\begin{theorem}{\em (\cite[Theorem 11]{ADX2017+})}
\label{thm:ec-pro}
If $G$ and $H$ are connected graphs, then we have $C_{\ecc}(G\cp H)=C_{\ecc}(G)+C_{\ecc}(H)-1$.
\end{theorem}

For any graph $G$, we denote by  $G^m$ be $m^{\rm th}$ power of $G$ with respect to Cartesian product, that is, the Cartesian product of $m$ copies of $G$. Now we have:

\begin{corollary}\label{co:ec-pro}
If $k\geq 1$, then $C_{\ecc}(G^{2^k})=2^kC_{\ecc}(G)-2^k+1$.
\end{corollary}

\proof By Theorem \ref{thm:ec-pro}, we have
\begin{eqnarray*}
C_{\ecc}(G^{2^k})
&=&2C_{\ecc}(G^{2^{k-1}})-1 \\
&=&2\Big[2C_{\ecc}(G^{2^{k-2}})-1\Big]-1 \\
&\vdots &\nonumber\\
&=&2^iC_{\ecc}(G^{2^{k-i}})-2^{i-1}-\cdots-2-1 \\
& \vdots & \\
&=&2^kC_{\ecc}(G)-2^k+1 \,.
\end{eqnarray*}
\qed

The Wiener complexity of Cartesian products is more involved. The distance formula~\eqref{eq:distance-formula} yields
\begin{eqnarray}
\tr_{G\cp H}((g,h))
&=&\sum\limits_{(g^{\prime},h^{\prime})\in V(G\cp H)}d_{G\cp H}((g,h),(g^{\prime},h^{\prime}))\nonumber \\
&=&\sum\limits_{(g^{\prime},h^{\prime})\in V(G\cp H)}(d_G(g,g^{\prime})+ d_H(h,h^{\prime}))\nonumber\\
&=&n(H)\sum\limits_{g^{\prime}\in V(G)}d_G(g,g^{\prime})+n(G)\sum\limits_{h^{\prime}\in V(H)}d_H(h,h^{\prime})\nonumber\\
&=&n(H)\tr_G(g)+n(G)\tr_H(h)\,,
\label{eq1}
\end{eqnarray}
a result deduced earlier in~\cite{AAKS2014}. Consequently,
\begin{equation}
\label{eq2}
C_W(G\cp H) = \big|\{n(H)\tr_G(g)+n(G)\tr_H(h):\ g\in V(G),h\in V(H)\}\big|\,,
\end{equation}
from which we immediately get:
\begin{equation}
\label{eq:upper-lower}
\max\{C_W(G),C_W(H)\} \leq C_W(G\cp H)\leq C_W(G)C_W(H)\,.
\end{equation}
The lower bound in~\eqref{eq:upper-lower} can be improved as follows.

\begin{proposition}
	\label{prp:lower-bound-Wiener}
	If $G$ and $H$ are graphs, then $$C_W(G \cp H) \geq C_W(G) + C_W(H) - 1\,.$$
\end{proposition}

\proof
Let $\tr(G) = \{ x_1, \ldots, x_s \}$ and $\tr(H) = \{ y_1, \ldots, y_t \}$, where $x_1 < \cdots < x_s$ and $y_1 < \cdots < y_t$. Then the set $X$ defined as $$\{x_1n(H) + y_1n(G), \ldots, x_1n(H) + y_tn(G), x_2n(H) + y_tn(G), \ldots, x_sn(H) + y_tn(G)\}$$ contains pairwise different integers, and $X \subseteq \tr(G \cp H)$ by~\eqref{eq1}. Since $$|X|  = s + t - 1 = C_W(G) + C_W(H) - 1,$$
the result follows.
\qed

We note in passing that Theorem~\ref{thm:ec-pro} and Proposition~\ref{prp:lower-bound-Wiener} immediately imply that if $C_W(G)\geq C_{\ecc}(G)$ and $C_W(H)\geq C_{\ecc}(H)$, then $C_W(G\cp H) \geq C_{\ecc}(G\cp H)$.

If $G$ is a graph and $H$ a graph with $C_W(H)=1$, then $C_W(G\cp H)=C_W(G)$, a result first reported in~\cite{AAKS2014}. Hence the lower bound in~\eqref{eq:upper-lower} is best possible, it coincides with that in Proposition~\ref{prp:lower-bound-Wiener}. On the other hand, the sharpness of the upper bound in~\eqref{eq:upper-lower} was not discussed in~\cite{AAKS2014}. To establish the sharpness also for the upper bound we introduce the following notion. A graph $G$ is {\em transmission indivisible} if $n(G)\nmid (\tr_G(u)-\tr_G(v))$ for every two distinct vertices $u,v\in V(G)$.

\begin{theorem}
\label{thm:coprime}
Let $G$ and $H$ be graphs. If at least one of $G$ and $H$ is transmission indivisible, and $\gcd(n(G),n(H))=1$, then $C_W(G\cp H)= C_W(G)C_W(H)$.
\end{theorem}

\proof
Let $\tr(G)=\{p_i:\ i\in [s]\}$ and $\tr(H)=\{q_j:\ j\in [t]\}$. Then in view of~\eqref{eq2},
$$C_W(G\cp H)=\big|\{n(H)p_i + n(G)q_j:\ i\in [s],j\in [t]\}\big|\,.$$
To prove the assertion of the theorem we need to show that $n(H)p_{i}+n(G)q_{j}\neq n(H)p_{i'}+n(G)q_{j'}$ for every $\{i,i'\}\in \binom{[s]}{2}$ and every $\{j,j'\}\in \binom{[t]}{2}$. Suppose on the contrary that for some such pairs $\{i,i'\}$ and $\{j,j'\}$ we have $n(H)p_{i}+n(G)q_{j} = n(H)p_{i'}+n(G)q_{j'}$, that is,
$n(H)(p_{i} - p_{i'}) = n(G)(q_{j'} - q_{j})$. Because $\gcd(n(G),n(H))=1$ we infer that $n(G)|(p_{i} - p_{i'})$ and $n(H)|(q_{j'} - q_{j})$. But this means that neither $G$ nor $H$ is transmission indivisible, a contradiction.
\qed

Note that if $G$ is transmission indivisible, then the transmissions of all the vertices of $G$ are pairwise different, that is, $G$ is transmission irregular. Although almost all graphs are not transmission irregular, an infinite family of transmission irregular trees was constructed in~\cite{AK2018} and an infinite family of transmission irregular trees of even order in~\cite{dobrynin-2019-b}. Moreover, an infinite family of transmission irregular 2-connected graphs was constructed in~\cite{dobrynin-2019} and an infinite family of 3-connected cubic transmission irregular graphs in~\cite{dobrynin-2019-c}.

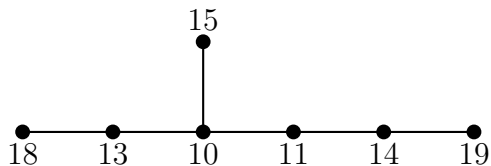
\begin{figure}[ht!]
\begin{center}
\begin{tikzpicture}[scale=0.8,style=thick]
\def\vr{3pt}
\begin{scope}
\path (2.5,0) coordinate (x1);
\path (4,0) coordinate (x2);
\path (5.5,0) coordinate (x3);
\path (7,0) coordinate (x4);
\path (8.5,0) coordinate (x5);
\path (10,0) coordinate (x6);
\path (5.5,1.5) coordinate (y3);

\draw (x1) -- (x2) -- (x3) -- (x4) -- (x5) --(x6);
\draw (x3) -- (y3) ;
\draw (x1)  [fill=black] circle (\vr);
\draw (x2)  [fill=black] circle (\vr);
\draw (x3)  [fill=black] circle (\vr);
\draw (x4)  [fill=black] circle (\vr);
\draw (x5)  [fill=black] circle (\vr);
\draw (x6)  [fill=black] circle (\vr);
\draw (y3)  [fill=black] circle (\vr);

\draw [below] (x1) node {$18$};
\draw [below] (x2) node {$13$};
\draw [below] (x3) node {$10$};
\draw [below] (x4) node {$11$};
\draw [below] (x5) node {$14$};
\draw [below] (x6) node {$19$};
\draw [above] (y3) node {$15$};
\end{scope}

\end{tikzpicture}
\end{center}
\caption{A transmission irregular but non-transmission indivisible tree}
\label{fig:tree}
\end{figure}

\begin{figure}[ht!]
\begin{center}
\begin{tikzpicture}[scale=0.8,style=thick]
\def\vr{3pt}
\begin{scope}
\path (0,3) coordinate (x9);
\path (2.5,1.8) coordinate (x11);
\path (3,0) coordinate (x10);
\path (6,1.5) coordinate (x13);
\path (5,3) coordinate (x8);
\path (0,4.5) coordinate (x15);
\path (2.5,4.5) coordinate (x12);
\path (5,4.5) coordinate (x14);
\draw (x10) -- (x9) -- (x11) -- (x8) -- (x13) --(x10);
\draw (x11) -- (x10) -- (x8) -- (x14);
\draw (x8) -- (x12) -- (x9) -- (x8);
\draw (x9) -- (x15);
\draw (x8)  [fill=black] circle (\vr);
\draw (x9)  [fill=black] circle (\vr);
\draw (x10)  [fill=black] circle (\vr);
\draw (x11)  [fill=black] circle (\vr);
\draw (x12)  [fill=black] circle (\vr);
\draw (x13)  [fill=black] circle (\vr);
\draw (x14)  [fill=black] circle (\vr);
\draw (x15)  [fill=black] circle (\vr);
\draw [right] (x8) node {$8$};
\draw [left] (x9) node {$9$};
\draw [left] (x10) node {$10$};
\draw [above] (x11) node {$11$};
\draw [above] (x12) node {$12$};
\draw [right] (x13) node {$13$};
\draw [right] (x14) node {$14$};
\draw [left] (x15) node {$15$};
\end{scope}

\begin{scope}[xshift = -8.5cm]
\path (0,3) coordinate (x9);
\path (3,0) coordinate (x10);
\path (6,1.5) coordinate (x13);
\path (5,3) coordinate (x8);
\path (0,4.5) coordinate (x15);
\path (2.5,4.5) coordinate (x12);
\path (5,4.5) coordinate (x14);
\draw (x10) -- (x9) -- (x8) -- (x13) --(x10);
\draw (x10) -- (x8) -- (x14);
\draw (x8) -- (x12) -- (x9);
\draw (x9) -- (x15);
\draw (x8)  [fill=black] circle (\vr);
\draw (x9)  [fill=black] circle (\vr);
\draw (x10)  [fill=black] circle (\vr);
\draw (x12)  [fill=black] circle (\vr);
\draw (x13)  [fill=black] circle (\vr);
\draw (x14)  [fill=black] circle (\vr);
\draw (x15)  [fill=black] circle (\vr);
\draw [right] (x8) node {$7$};
\draw [left] (x9) node {$8$};
\draw [left] (x10) node {$9$};
\draw [above] (x12) node {$10$};
\draw [right] (x13) node {$11$};
\draw [right] (x14) node {$12$};
\draw [left] (x15) node {$13$};
\end{scope}

\end{tikzpicture}
\end{center}
\caption{Two of the smallest interval irregular graphs}
\label{fig:sporadic}
\end{figure}
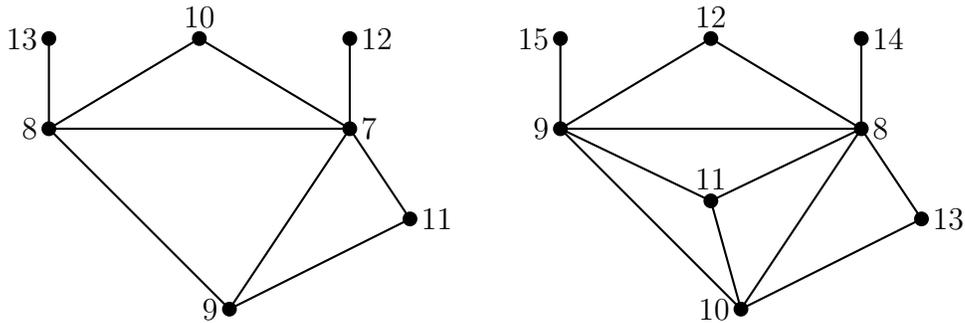

A transmission irregular graph need not be transmission indivisible. In Fig.~\ref{fig:tree} a transmission irregular but non-transmission indivisible tree of order $7$ is shown, where the transmission is given for each vertex. Sporadic transmission indivisible graphs of order $7$ and $8$ are shown in Fig.~\ref{fig:sporadic}, where along with each vertex its transmission is stated.

Note that the transmission of the vertices in each of the two examples from Fig.~\ref{fig:sporadic} are consecutive integers from the intervals $[7..13]$ and $[8..15]$, respectively, which makes these examples particularly interesting. Such graphs were named \textit{interval (transmission) irregular graphs} in~\cite{dobrynin-2019}. Interval irregular graphs are transmission indivisible and hence Theorem~\ref{thm:coprime} applies. We have checked by computer that there are 1, 2, 13, and 0 interval irregular graphs on 7, 8, 9, and 10 vertices, respectively. In addition, Dobrynin~\cite{dobrynin-2019} reports that {\s there exist at least 207 interval irregular $2$-connected graphs of order 11}. Their respective intervals of transmissions are $[13..23]$ (154 graphs), $[15..25]$ (51 graphs), and $[17..27]$ (2 graphs). The existence of an infinite family of interval irregular graphs is an open problem.

With interval irregular graphs in hands (and hence with transmission indivisible ones) the following result makes sense.

\begin{corollary}
\label{cor:primes}
If $G$ is a transmission indivisible graph, then there exits a family of graphs $\{H_i\}_{i\ge 1}$  such that $C_W(G\cp H_i)= C_W(G)C_W(H_i)$.
\end{corollary}

\proof
Let $\{p_i\}_{i\ge 1}$ be a set of primes each larger than $n(G)$ and let $H_i$ be a graph of order $p_i$. Then $\gcd(n(G),n(H_i))=1$ and the result follows from Theorem~\ref{thm:coprime}.
\qed

\begin{figure}[ht!]
	\begin{center}
		\begin{tikzpicture}[scale=0.8,style=thick]
		\def\vr{3pt}
		\begin{scope}
		\path (0,0) coordinate (x8);
		\path (1.5,0) coordinate (x10);
		\path (3,0) coordinate (x9);
		\path (4.5,0) coordinate (x7);
		
		\path (0,1.5) coordinate (x4);
		\path (1.5,1.5) coordinate (x0);
		\path (3,1.5) coordinate (x6);
		\path (4.5,1.5) coordinate (x2);
		
		\path (0,-1.5) coordinate (x1);
		\path (1.5,-1.5) coordinate (x5);
		\path (4.5,-1.5) coordinate (x3);

		\draw (x9) -- (x10) -- (x5) -- (x1) -- (x8) -- (x0) -- (x10) -- (x6) --(x9) -- (x2) -- (x7) -- (x3) -- (x9) -- (x7);
		\draw (x4) -- (x10) -- (x8) -- (x4) -- (x0) -- (x6) --(x2);

		\draw (x0)  [fill=black] circle (\vr);
		\draw (x1)  [fill=black] circle (\vr);
		\draw (x2)  [fill=black] circle (\vr);
		\draw (x3)  [fill=black] circle (\vr);
		\draw (x4)  [fill=black] circle (\vr);
		\draw (x5)  [fill=black] circle (\vr);
		\draw (x6)  [fill=black] circle (\vr);
		\draw (x7)  [fill=black] circle (\vr);
		\draw (x8)  [fill=black] circle (\vr);
		\draw (x9)  [fill=black] circle (\vr);
		\draw (x10)  [fill=black] circle (\vr);

		\draw [above] (x0) node {$18$};
		\draw [below] (x1) node {$26$};
		\draw [above] (x2) node {$22$};
		\draw [below] (x3) node {$24$};
		\draw [above] (x4) node {$20$};
		\draw [below] (x5) node {$21$};
		\draw [above] (x6) node {$17$};
		\draw [right] (x7) node {$23$};
		\draw [left] (x8) node {$19$};
		\draw [below] (x9) node {$16\quad$};
		\draw [below right] (x10) node {$14$};
		\end{scope}
		
		\end{tikzpicture}
	\end{center}
	\caption{A $2$-connected transmission indivisible graph which is not interval irregular}
	\label{fig:indivisible}
\end{figure}
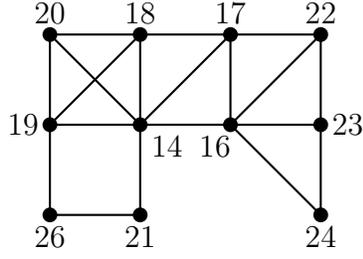

By a computer search (using~\cite{nauty}) we have checked that the class of transmission indivisible graphs is strictly larger than the class of interval irregular graphs. There are no such examples on up to and including $10$ vertices.  However, a bit surprisingly, there are 221 graphs on $11$ vertices that are transmission indivisible but not interval irregular. Among them there are no trees, but one finds 14 graphs which are $2$-connected, an example can be seen in Fig.~\ref{fig:indivisible}.

A graph $G$ is \emph{arithmetic transmission} if  the ordered elements of $\tr(G)$ form an arithmetic progression. Moreover, if $\tr(G)$ has step $a$, we say that $G$ has step $a$. (In Subsection~\ref{subsec:larger} see an example of arithmetic transmission graph with step $4$.) Below we present a result in which the lower bound is attained in Proposition~\ref{prp:lower-bound-Wiener}.

\begin{proposition}
	\label{ttap}
	Let $G$ and $H$ be arithmetic transmission graphs. If $\tr(H)\subseteq \tr(G)$ and $n(G) = n(H)$, then $C_W(G\cp H)=C_W(G)+C_W(H)-1$.\end{proposition}
\proof Set $n = n(G) = n(H)$ and $\tr(G)=\{x_1,\ldots,x_k\}$. We may without loss of generality assume that $\tr(H)=\{x_1,\ldots,x_j\}$, where $j \leq k$. Then $|\tr(G)| = k = \frac{x_k - x_1}{a} + 1$ and $|\tr(H)| = j = \frac{x_j - x_1}{a} + 1$.  By \eqref{eq2}, $\tr(G \cp H) = \{ n (x_1 + x_1), \ldots, n (x_k + x_j) \}$ and hence
\begin{align*}
|\tr(G\cp H)| & = \frac{(x_k + x_j) - (x_1 + x_1)}{a} + 1 \\
              & = \frac{x_k - x_1}{a} + \frac{x_j - x_1}{a} + 1 + 1 - 1 \\
              & = |\tr(G)| + |\tr(H)| - 1\,.
\end{align*}

\vspace{-21pt}
\qed

 Taking $H = G$ in Proposition~\ref{ttap} and using a similar technique as that in the proof of Corollary \ref{co:ec-pro}, we have the following result.
\begin{corollary}
	\label{co:wc-pro}
	Let $G$ be an arithmetic transmission graph and $k\geq 0$ be an integer. Then $C_{W}(G^{2^k})=2^kC_{W}(G)-2^k+1$.
\end{corollary}

\vspace{10pt}

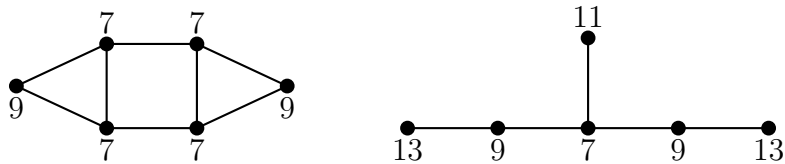
\begin{figure}[ht!]
\begin{center}
\begin{tikzpicture}[scale=0.8,style=thick]
\def\vr{3pt}
\begin{scope}
\path (-4,0.7) coordinate (u);
\path (-2.5,0) coordinate (w1);
\path (-2.5,1.4) coordinate (w2);
\path (-1,0) coordinate (w3);
\path (-1,1.4) coordinate (w4);
\path (0.5,0.7) coordinate (v);

\path (2.5,0) coordinate (x1);
\path (4,0) coordinate (x2);
\path (5.5,0) coordinate (x3);
\path (7,0) coordinate (x4);
\path (8.5,0) coordinate (x5);
\path (5.5,1.5) coordinate (y3);

\draw (u) -- (w1) -- (w3) -- (v);
\draw (u) -- (w2) -- (w4) -- (v);
\draw  (w1) -- (w2);
\draw  (w3) -- (w4);
\draw (x1) -- (x2) -- (x3) -- (x4) -- (x5);
\draw (x3) -- (y3) ;
\draw (x1)  [fill=black] circle (\vr);
\draw (x2)  [fill=black] circle (\vr);
\draw (x3)  [fill=black] circle (\vr);
\draw (x4)  [fill=black] circle (\vr);
\draw (x5)  [fill=black] circle (\vr);
\draw (y3)  [fill=black] circle (\vr);

\draw (u)  [fill=black] circle (\vr);
\draw (w1)  [fill=black] circle (\vr);
\draw (w2)  [fill=black] circle (\vr);
\draw (w3)  [fill=black] circle (\vr);
\draw (w4)  [fill=black] circle (\vr);
\draw (v)  [fill=black] circle (\vr);

\draw [below] (x1) node {$13$};
\draw [below] (x2) node {$9$};
\draw [below] (x3) node {$7$};
\draw [below] (x4) node {$9$};
\draw [below] (x5) node {$13$};
\draw [above] (y3) node {$11$};

\draw [below] (u) node {$9$};
\draw [below] (w1) node {$7$};
\draw [above] (w2) node {$7$};
\draw [below] (w3) node {$7$};
\draw [above] (w4) node {$7$};
\draw [below] (v) node {$9$};
\end{scope}

\end{tikzpicture}
\end{center}
\caption{ Graphs $G$ (right) and $H$ (left) with ${\rm Tr(H)}\subset \tr(G)$}
\label{fig:ttap-graph}
\end{figure}

In Fig.~\ref{fig:ttap-graph} graphs $G$ and $H$ are shown which satisfy the conditions of Proposition~\ref{ttap}.

Clearly, any interval irregular graph is arithmetic transmission as a special case. Hence  Proposition \ref{ttap} and Corollary \ref{co:wc-pro} also apply for interval irregular graphs. Next we provide a result on the Cartesian product of arithmetic transmission graphs for which the upper bound in~\eqref{eq:upper-lower} is attained.

\begin{theorem}
	\label{two-tap}
	Let $G$ and $H$ be arithmetic transmission graphs with steps $a$ and $b$, respectively. If $n(H)(C_W(G)-1)a<n(G)b$, then $C_W(G\cp H)=C_W(G)C_W(H)$. \end{theorem}
\proof Let $\tr(G)=\{x_1,\ldots,x_k\}$ and ${\rm Tr(H)}=\{y_1,\ldots,y_j\}$. Thus $C_W(G)=k$, $C_W(H)=j$ and
\begin{equation}
C_W(G\cp H) = \big|\{n(H)x_p+n(G)y_q:\ p\in [k],q\in [j]\}\big|\nonumber
\end{equation} by~\eqref{eq2}. Let $A_i=\{n(H)x_1+n(G)y_i,\ldots,n(H)x_k+n(G)y_i\}$ for $i\in [j]$. Therefore, we have $\tr(G\cp H)=\bigcup\limits_{i=1}^{j}A_i$. From the assumption, we have $n(H)(k-1)a<n(G)b$. It follows that
\begin{align*}
n(H)x_k+n(G)y_i & = n(H) (x_1 + (k-1) a) + n(G) (y_{i+1} - b) \\
                & = (n(H) x_1 + n(G) y_{i+1}) + (n(H) (k-1) a - n(G) b)\\
                & <  n(H)x_1+n(G)y_{i+1}\,.
\end{align*}
This implies that $\max\{s:s\in A_i\}<\min\{s:s\in A_{i+1}\}$, hence $A_s\cap A_t=\emptyset$ for any $s,t\in [k]$, $s \neq t$. We conclude that $C_W(G\cp H)=kj$.\qed

\section{Comparing $C_{\ecc}(G)$ with $C_W(G)$}
\label{sec:compare}

It was proved in~\cite{ADX2017+} that $C_{\ecc}(G)\leq\lceil\frac{n(G)}{2}\rceil$ for any graph $G$. Therefore, if $G$ is a transmission irregular graph, then $C_W(G)>C_{\ecc}(G)$. Actually, this is a phenomena that is very common as the next result shows.

\begin{proposition}
\label{prp:for-almosst-all}
For almost all graphs $G$, we have $C_{\ecc}(G) \le C_W(G)$.
\end{proposition}

\proof
It is well-known that almost all graphs have diameter $2$. So let $G$ be a graph with $\diam(G) = 2$. Then $C_{\ecc}(G) \le 2$. There is nothing to show if $C_{\ecc}(G) = 1$, so let $C_{\ecc}(G) = 2$ in which case we have ${\rm Ec}(G) = \{1,2\}$. But then $G$ contains at least one vertex of degree $n(G) - 1$, and at least one vertex of smaller degree. Since their transmissions are different by Proposition~\ref{prp:2-ec}, $C_W(G)\ge 2$.
\qed

We next show that the difference $C_W(G) - C_{\ecc}(G)$ can be arbitrarily large. For this recall that the $d$-cube $Q_d$, $d \geq 1$, has the vertex set $\{0, 1\}^d$, two vertices in $Q_d$ are adjacent if they differ in precisely one coordinate. Let $Q_d^-$ be the graph obtained from $Q_d$ by removing an arbitrary vertex.

\begin{proposition}
	\label{prp:cube}
	If $d \geq 2$, then $C_W(Q_d^-) - C_{\ecc}(Q_d^-) = d-2$.
\end{proposition}

\proof
We may without loss of generality assume that $V(Q_d^-) = V(Q_d) \setminus\{0^d\}$. Then $\ecc_{Q_d^-}(1^d) = d-1$ and $\ecc_{Q_d^-}(v) = d$ for all other vertices $v$ of $Q_d^-$. Thus $C_{\ecc}(Q_d^-) = 2$.

Let $x \in V(Q_d^-)$. Then $\tr_{Q_d^-}(x) = \tr_{Q_d}(x) \setminus\{d_{Q_d}(x, 0^d)\}$. Since there are precisely $d$ different values of $d_{Q_d}(x, 0^d)$, we conclude that $\tr(Q_d^-) = d$.
\qed

\begin{corollary}
	\label{cor:larger-diameter}
	If $d \geq 2$, then there exists a graph $G$ with $\diam(G) = d$ such that $C_W(G) > C_{\ecc}(G)$.
\end{corollary}

\proof
If $d = 2$, then consider the paw graph $P$ (that is, the graph obtained from a triangle by attaching a leaf to one of its vertices) for which ${\rm Ec}(P) = \{ 1,2 \}$ and $\tr(P) = \{3,4,5\}$ hold. For $d \geq 3$, apply Proposition~\ref{prp:cube}.
\qed

\subsection{ More graphs $G$ with  $C_W(G)\geq C_{\ecc}(G)$}
\label{subsec:equality-case}

Let $G$ be a graph. Then, by definition, $C_{\ecc}(G)=C_W(G)=1$ if and only if $G$ is a self-centered, transmission regular graph. In the next result we present several additional classes of graphs $G$ for which $C_W(G)=C_{\ecc}(G)$ holds. To state the result, we need some further definitions. If the eccentricity of the vertices of a self-centered graph $G$ is $k$, we say that $G$ is {\em $k$-self-centered}. A graph $G$ is {\em bidegreed} if all vertices of $G$ have one of two possible degrees, cf.~\cite{myrvold-1987}. (If $x$ is a vertex of a bidegreed graph $G$, then $\deg_G(x) \in\{\delta(G), \Delta(G)\}$.) Let finally $G^{(k*)}$, $k\ge 1$, denote the graph obtained from $G$ by attaching $k$ pendant vertices to each vertex of $G$.

\begin{proposition}
\label{prop:W=Ec-classes}
\begin{enumerate}
\item[(i)] If $G$ is a regular, $2$-self-centered graph, then $C_W(G)=C_{\ecc}(G)=1$.
\item[(ii)] If $G$ {\s is} a bidegreed, non-self-centered graph with $\diam(G) = 2$, then $C_W(G)=C_{\ecc}(G) = 2$.
\item[(iii)] If $G$ is a regular or bidegreed graph obtained from $K_n$ by removing $k$ edges, $k\in [\lfloor\frac{n}{2}\rfloor]$, then $C_W(G)=C_{\ecc}(G)$.
\item[(iv)] If $G$ is a vertex-transitive graph and $k\ge 1$, then $C_{\ecc}(G^{(k*)}) = C_{W}(G^{(k*)}) = 2$.
\end{enumerate}
\end{proposition}

\proof
(i) As $G$ is self-centered, $C_{\ecc}(G)=1$. Combining  Proposition~\ref{prp:2-ec} and the fact that $G$ is regular, we have $C_W(G)=1$.

(ii) Since $G$ is not self-centered and $\diam(G) = 2$, we have $C_{\ecc}(G)=2$ and at least one vertex must have eccentricity $1$, that is, $\Delta(G)=n(G)-1$. As $G$ is bidegreed, all the vertices that have degree smaller than $n(G)-1$ must have the same degree. In view of Proposition~\ref{prp:2-ec}, we have $C_W(G)=2$.

(iii) If $G$ is regular, then since $k\le \lfloor\frac{n}{2}\rfloor$, $G$ is obtained from $K_n$ be removing a perfect matching (in which case $n$ is even and $k = n/2$). Then the assertion follows by~(i). Otherwise $G$ is not regular. But then $G$ is bidegreed and hence $G$ fulfils the assumption of~(ii).

(iv) Since $G$ is vertex-transitive, $G$ is a transmission regular, self-centered graph. Assume that $\diam(G)=d$, $n = n(G)$, and let $V(G)=\{v_1,\ldots,v_n\}$. Let $V(G^{(k*)})=V(G)\cup\{v_i^{(j)}:\ i\in [n],j\in [k]\}$ with $v_iv_i^{(j)}\in E(G^{(k*)})$ for $j\in [k]$. From the structure of $G^{(k*)}$ we have $\ecc_{G^{(k*)}}(v_i)=d+1$ and $\ecc_{G^{(k*)}}(v_i^{(j)})=d+2$ for every $i\in [n]$, $j\in [k]$. Hence $C_{\ecc}(G^{(k*)})=2$.

On the other hand, $\tr_{G^{(k*)}}(v_i)$ is the same for all $i\in [n]$. Since $v_i^{(j)}$ is a pendant vertex, we see that $\tr_{G^{(k*)}}(v_i^{(j)})=\tr_{G^{(k*)}}(v_i)+(n+1)k-2$ is the same for every $j\in [k]$. Thus $C_{W}(G^{(k*)})=2$.
\qed

The class of graphs from Proposition~\ref{prop:W=Ec-classes} (i) contains vertex-transitive graphs as a proper subclass. For instance, if $G$ is an arbitrary regular graph that is not vertex-transitive, then the join of two copies of $G$ is a regular, $2$-self-centered graph, but not vertex-transitive. (The {\em join} of graphs $G$ and $H$ is obtained from the disjoint union of $G$ and $H$ by adding all possible edges between vertices of $G$ and vertices of $H$.)

To show that there exist graphs that have the same Wiener complexity and eccentric complexity which is arbitrary large, Cartesian and lexicographic product graphs can be used. We have already defined the Cartesian product. The {\em lexicographic product} $G\circ H$ of graphs $G$ and $H$ also has the vertex set $V(G)\times V(H)$, vertices $(g,h)$ and $(g',h')$ being adjacent if either $gg'\in E(G)$, or $g=g'$ and $hh'\in E(H)$.

\begin{theorem}
\label{thx-lex-and-Cart}
(i) If $G$ is a graph with $C_W(G)=C_{\ecc}(G)$, and $H$ is a self-centered, transmission regular graph, then
$$C_{\ecc}(G\cp H)=C_{W}(G\cp H) = C_W(G)\,.$$
(ii) If $H$ is a regular graph and $n\ge 4$, then
$$C_W(P_n\circ H) = C_{\ecc}(P_n\circ H) = \left\lceil \frac{n}{2}\right\rceil\,.$$
\end{theorem}

\proof
(i) The assumption that $H$ is a self-centered, transmission regular graph, means that $C_{\ecc}(H)=C_W(H)=1$. Then Theorem~\ref{thm:ec-pro} implies that
$$C_{\ecc}(G\cp H) = C_{\ecc}(G) + C_{\ecc}(H) - 1 = C_{\ecc}(G) = C_{W}(G)\,.$$
On the other hand, from~\eqref{eq:upper-lower} we get $C_{W}(G\cp H) = C_W(G)$ because $C_W(H)=1$.

(ii)  Let $V(P_n) = \{v_1, \ldots, v_n\}$ with natural adjacency relation.

 Consider first the case when $H=K_1$. Then $P_n\circ H = P_n\circ K_1 = P_n$. Clearly, $\ecc_{P_n}(v_i)=|i|_n$, where $|i|_n=\max\{i-1,n-i\}$. Consequently $C_{\ecc}(P_n) = \left\lceil \frac{n}{2}\right\rceil$. Moreover, if $i\in [n]$, then $\tr_{P_n}(v_i)={i \choose 2}+{n-i+1 \choose 2}$. In particular, if $i,j\le  \left\lceil \frac{n}{2}\right\rceil$, $i\ne j$, then $\tr_{P_n}(v_i) \ne \tr_{P_n}(v_j)$, and $\tr_{P_n}(v_i) = \tr_{P_n}(v_{n-i+1})$. It follows that $C_{W}(P_n) = \left\lceil \frac{n}{2}\right\rceil$.

Let now $H$ be an arbitrary regular graph and consider the lexicographic product $P_n\circ H$. Note first that $d_{P_n\circ H}((v_i,h), (v_i,h'))\le 2$ for any vertices $h,h'\in V(H)$. Moreover, $d_{P_n\circ H}((v_i,h), (v_j,h')) = d_{P_n}(v_i,v_j)$ for $i\ne j$. Since $n\ge 4$ it follows that $\ecc_{P_n\circ H}((v_i,h)) = \ecc_{P_n}(v_i) = |i|_n$ and hence $C_{\ecc}(P_n\circ H) = \left\lceil \frac{n}{2}\right\rceil$.

Consider now vertices $(v_i,h)$ and $(v_i,h')$ of $P_n\circ H$ for some $i\in [n]$ and $h,h'\in V(H)$, $h\ne h'$. Let $V_i = \{(v_i,x):\ x\in V(H)\}$. Since $H$ is regular, both vertices $(v_i,h)$ and $(v_i,h')$ have the same number of neighbors in $V_i$. Moreover, the distance between them and their non-neighbors in $V_i$ is $2$. Since we already observed that $d_{P_n\circ H}((v_i,h), (v_j,h'')) = d_{P_n}(v_i,v_j) = d_{P_n\circ H}((v_i,h'), (v_j,h''))$ for every $h'' \in V(H)$, we get that $\tr_{P_n\circ H}((v_i,h)) = \tr_{P_n\circ H}((v_i,h'))$, that is, the vertices of $V_i$ have the same transmission. Moreover, by the argument from the case $H=K_1$ we also get that if $i,j\le  \left\lceil \frac{n}{2}\right\rceil$, $i\ne j$, and $h\in V(H)$, then $\tr_{P_n\circ H}((v_i,h)) \ne \tr_{P_n\circ H}((v_j,h))$, and $\tr_{P_n\circ H}((v_i,h)) = \tr_{P_n\circ H}((v_{n-i+1}, h))$. We conclude that $C_{W}(P_n\circ H) = \left\lceil \frac{n}{2}\right\rceil$.
\qed

For a connected graph $G$, the set of vertices at the given distance from $C(G)$ is called a \textit{distance-level} of $G$.  The set $L_i$ of vertices at distance $i$ from $C(G)$ is called \textit{$i$-distance-level} of $G$ for $i\in [\rad(G)]$. A tree $T$ is \emph{center-regular} if the vertices in  every distance-level have the same degree. Note that this in particular implies that if $T$ is bicentered, then the central vertices have the same degree.

\begin{proposition}
	\label{prp:center-regular-tree}
	If $T$ is a center-regular tree, then $C_W (T) = C_{\ecc}(T)$.
\end{proposition}

\proof
Let $k = \rad(T)$ and consider the $i$-distance-levels $L_i$ of $T$, $0 \leq i \leq \rad(T)$. Let $u$ and $v$ be arbitrary vertices from  $L_i$. Then since $T$ is center-regular, there exists an automorphism $\varphi \in {\rm Aut}(T)$ such that $\varphi(u) =  v$. This implies that all the vertices of $L_i$ have the same eccentricity as well as the same transmission. Hence, $C_{\ecc} (T) \leq k+1$ and $C_{W}(T) \leq k+1$. On the other hand, it is obvious that $C_{\ecc} (T) \geq k+1$ and $C_{W}(T) \geq k+1$.
\qed

The family of center-regular trees includes as a special case the recently introduced \emph{degree-eccentricity regular (DE-regular for short)} trees~\cite{XGG2018}. The degree-eccentricity regular trees are defined as the trees $T$ in which $\deg_T(v)+ \ecc_T(v)$ is a fixed constant for every vertex $v\in V(T)$. Clearly, such a tree is center-regular. In~\cite{XGG2018}, among other results, all the molecular $DE$-regular trees were completely characterized.

\begin{remark}
	The result from Proposition~\ref{prp:center-regular-tree} can be generalized as follows. Let $G$ be a graph such that for every two vertices in a distance-level with respect to $C(G)$, there exists an automorphism mapping one vertex to the other. If, in addition, transmissions are pairwise different for distance-levels, then $C_W(G) = C_{\ecc}(G)$.
\end{remark}
 \begin{theorem}
 If $T$ is a tree, then $C_{\ecc}(T) \leq C_W(T)$. Moreover, equality holds if and only if $T$ is center-regular tree.
\end{theorem}

\proof  Set $n = n(T) \ge 2$ and consider the following two cases.

\medskip\noindent
{\bf Case 1.} $|C(T)| = 1$. \\
In this case $\diam(T) = 2\,\rad(T)$ holds. Therefore, in view of~\eqref{eq:C_ec}, we need to prove that the number of distinct values of $\tr(v)$ is greater than or equal to $C_{\ecc}(T) = \rad(T) + 1$. Let $C(T) = \{r\}$ and consider $T$ as a tree rooted in $r$. Let $c(v)$ be the order of the subtree rooted at the vertex $v$ and containing all the vertices $x$ such that $v$ lies on the shortest $r,x$-path. Note that $v$ itself lies in this tree. For example $c(r) = n$ and if $v$ is a leaf, then $c(v) = 1$.

By definition, there are at least two disjoint paths of length $\rad(T)$ starting at $r$. Consider such a path
$$P:v_1 \rightarrow v_2 \rightarrow \cdots \rightarrow v_{\rad(T)} \rightarrow r\,,$$
where $v_1$ is a leaf, and $c(v_{\rad(T)}) < n / 2$ holds. As there are at least two radial paths starting from $r$, such a path always exists.

Because the distances from the vertices below $v_i$ increase by $1$ and all others decrease by $1$, we infer that for every consecutive vertices $v_{i+1}$ and $v_{i}$ of $P$ it holds
$$\tr(v_{i+1}) = \tr(v_i) + 2 c(v_i) - n~~for~~i\in [\rad(T)-1]\,,$$
and that $\tr(r) = \tr(v_{\rad(T)}) + 2 c(v_{\rad(T)}) - n$ holds. Since $c(v_i) < n / 2$ holds for every $i \in [\rad(T)]$, we have strict chain of inequalities
$$\tr(v_1) > \tr(v_2) > \cdots > \tr(v_{\rad(T)}) > \tr(r)\,.$$
This already implies that these are  at least $\rad(T) + 1$ distinct values of transmissions---which we wanted to show.

If $T$ is center-regular tree, then the equality holds by Proposition \ref{prp:center-regular-tree}. Conversely, suppose that the equality holds for a tree $T$ and consider again $T$ rooted in its center $r$. Then for every two radial paths starting from $r$ there is a an automorphism mapping one path onto the other. We can show by simple induction that all vertices on the same distance from the root need to have the same degree. Namely, all the leaves on $\rad(T)$ level have the same degree 1. On the level $\rad(T) - 1$, as the numbers $c(v)$ need to be equal for all paths, it directly follows that these vertices have the same degree. We can continue this until we reach the root $r$.

\medskip\noindent
{\bf Case 2.} $|C(T)| = 2$. \\
In this case we need to prove (again in view of~\eqref{eq:C_ec}) that the number of distinct values of $\tr(v)$ is greater than or equal to $\rad(T)$.
We can use a  parallel technique as in Case $1$ for the  two center vertices, and also conclude that the equality holds if and only if the tree $T$ is center-regular.
\qed

\subsection{Graphs $G$ with $C_{\ecc}(G) > C_W(G) $}
\label{subsec:larger}

In this section, we construct graphs $G$ in which the Wiener complexity is arbitrarily smaller than the eccentric complexity. The existence of such graphs is not obvious at the first sight. Consider the following example.

 Let $Z_k$, $k\geq 1$, be the graph obtained by attaching a pendant vertex to each of two diametrical vertices in a cycle $C_{2k+2}$. Let $u$ and $v$ be the degree $3$ vertices of $Z_k$, let $u'$ and $v'$ be its respective neighbors of degree $1$, and let $x_1, \ldots, x_k$ and $y_1, \ldots, y_k$ be the $u, v$-paths in $Z_k$. Since the transmission of a vertex in $C_n$ is $\lfloor\frac{n^2}{4}\rfloor$, see~\cite{Plesn1984}, we have
\begin{align*}
	\label{eq:C_W-for-G-by-G}
	\tr_{Z_k}(x_i) & = \tr_{Z_k}(y_i) \\
	&=\tr_{C_{2k+2}}(x_i)+(i+1)+(k+1-i+1)  \\
	& = (k+1)^2+k+3,
\end{align*}
for $i\in [k]$. Moreover, $\tr_{Z_k}(u)=\tr_{Z_k}(v)=(k+1)^2+k+3$, and $\tr_{Z_k}(u')=\tr_{Z_k}(v')=(k+1)^2+3k+5$. On the other hand, ${\rm Ec}(Z_k)=\{k+1,k+2,k+3\}$. Thus $C_{\ecc}(Z_k)-C_W(Z_k)=1$.

 Each of the graphs $Z_k$ leads to another infinite family of graphs for which the eccentricity complexity exceeds the Wiener complexity. For this sake we  recall the following result.

\begin{corollary}{\em (\cite[Corollary 3.2]{AAKS2014})}
	\label{cor:andova}
	If $G$ is a graph and $H$ a graph with $C_W(H) = 1$, then $C_W(G \cp H) = C_W(G)$.
\end{corollary}

\begin{proposition}
	\label{prop:consecutive-infinite}
	If $G$ is a graph with $C_{\ecc}(G) > C_W(G)$ and $d \geq 1$, then $$C_{\ecc}(G \cp Q_d) > C_W(G \cp Q_d)\,.$$
\end{proposition}

\proof
By Corollary~\ref{cor:andova}, $C_W(G \cp Q_d) = C_W(G)$. On the other hand, Theorem~\ref{thm:ec-pro} implies that $C_{\ecc}(G\cp Q_d) = C_{\ecc}(G) + C_{\ecc}(Q_d) - 1 = C_{\ecc}(G)$.
\qed

Note that in the proof of Proposition~\ref{prop:consecutive-infinite}, the family of hypercubes could be replaced by an arbitrary family of graphs $\{H_i\}_{i \geq 1}$ with $C_W(H_i) = 1$.

We have thus seen that there are infinitely many graphs with eccentric complexity larger than the Wiener complexity. In the above families, this difference was $1$. We now demonstrate that the difference can be arbitrarily large. Let $Z=Z_1$ and recall that $Z^m$ is $m^{\rm th}$ power of $Z$ with respect to Cartesian product.

\begin{proposition}\label{Z-power}
If $k \geq 0$, then
$C_{\ecc}(Z^{2^k}) - C_W(Z^{2^k}) = 2^{k}$.
\end{proposition}

\proof
As we have observed above, ${\rm Ec}(Z) = \{2, 3, 4\}$ and $\tr(Z) = \{ 8, 12 \}$  with $C_{\ecc}(Z)=3$ and $C_W(Z)=2$. By Corollary \ref{co:ec-pro}, we have $C_{\ecc}(Z^{2^k})=3 \cdot 2^k-2^k+1=2^{k+1}+1$. From Corollary \ref{co:wc-pro}, we have $C_W(Z^{2^k})=2\cdot 2^k-2^k+1=2^k+1$.
The conclusion now follows immediately.
\qed

Based on   Corollaries~\ref{co:ec-pro} and~\ref{co:wc-pro}, Proposition~\ref{Z-power} can be generalized as follows.
\begin{corollary} Let $G$ be an arithmetic transmission graph with $C_{\ecc}(G)> C_W(G)$  and $k\geq 0$ be an integer. Then $C_{\ecc}(G^{2^k})> C_W(G^{2^k})$.\end{corollary}
Using the same argument as in the proof of Proposition~\ref{prop:consecutive-infinite}, we infer that for any positive integer $N$ there exists an infinite family of graphs $\{ H_i \}_{i \geq 1}$ such that $C_{\ecc}(H_i) - C_W(H_i) > N$.

As introduced in \cite{XDKL2020}, a graph $G$ is \textit{universally diametrical} (UD for short) if  there exist diametrical vertices $u$ and $v$ of $G$ such that $\text{Ecc}_G(w)\cap\{u,v\}\neq\emptyset$ for every vertex $w\in V(G)\setminus \{u,v\}$, that is, at least one of $u$ and $v$ is eccentric to $w$. Here $\text{Ecc}_G(w) = \{ u \in V(G) :\ d_G(u,w) = \ecc_G(w) \}$ is the \textit{eccentric set} (\cite{XLDK2018}) of $w$ in $G$. Here the vertices $u$ and $v$ form a \textit{universally diametrical pair} in $G$. A universally diametrical graph $G$ with a universally diametrical pair $u,v$ is called a \textit{$k$-$(u,v)$-universally diametrical} (or $k$-$(u,v)$-UD for simplicity) if $d_G(u,v)=\diam(G)=k$. Since any tree is an UD-graph, UD-graphs can be viewed as a generalization of trees. Moreover, the graph $Z_k$ defined as above is also an UD-graph with two universally diametrical vertices having equal transmissions such that $C_{\ecc}(Z_k)>C_W(Z_k)$. Next we present a method for constructing new graphs with $C_{\ecc}>C_W$ from UD-graphs.

\begin{theorem}\label{New-1}
Let $G_0$ be a $(v,v')$-UD graph with $\tr_G(v)=\tr_G(v^{\prime})$, let $T$ be a center-regular tree with $C(T) = \{x\}$, and let $T'\cong T$ with $C(T') = \{x'\}$. Let $G$ be a graph obtained from $G_0$, $T$, and $T'$ by identifying the vertices $v$ and $x$ and identifying the vertices $v'$ and $x'$. If $C_{\ecc}(G_0)>C_{W}(G_0)$, then $C_{\ecc}(G)>C_{W}(G)$.
\end{theorem}

\proof
Set $k = \diam(G_0)$, $r = \ecc_T(x) = \rad(T)$, $n_0 = n(G_0)$, and $n_1 = n(T) (= n(T'))$. For convenience, we still denote by $v$ and $v^{\prime}$ the  vertices of $G$ obtained by identifying $v$ with $x$ and by identifying $v^{\prime}$ with $x^{\prime}$, respectively. From the structure of $G$,  for any vertex $y\neq v$ from the $i$-distance-level of $T$ for $i\in [r]$, we have $\ecc_G(y)=r+k+i$ and  \begin{align*}
\tr_G(y) & = \sum\limits_{z\in V(T)}d_G(y,z)+\sum\limits_{z\in V(G_0)\setminus\{v,v^{\prime}\}}d_G(y,z)+\sum\limits_{z\in V(T^{\prime})}d_G(y,z) \\
& = \tr_T(y)+\sum\limits_{z\in V(G_0)\setminus\{v,v^{\prime}\}}\Big[d_T(y,v)+d_{G_0}(v,z)\Big] \\
&\quad + \sum\limits_{z\in V(T)}\Big[d_T(y,v)+k+d_{T^{\prime}}(v^{\prime},z)\Big] \\
&= \tr_T(y)+(n_0-2)i+\tr_{G_0}(v)-k+n_1(i+k)+\tr_{T^{\prime}}(v^{\prime})\\
&=\tr_T(y)+\tr_T(v)+\tr_{G_0}(v)+(n_1+n_0-2)i+(n_1-1)k\,.
\end{align*}
Therefore the transmission and the eccentricity can be uniquely determined by the value of $i$ for all vertices in the $i$-distance-level of $T$ in $G$. The same applies to the vertices from the $i$-distance-level of $T^{\prime}$ in $G$. From $\tr_G(v)=\tr_G(v^{\prime})$, we have
$$|\{\tr_G(y):y\in V(T)\cup V(T^{\prime})\setminus\{v,v^{\prime}\}\}|=|\{\ecc_G(y):y\in V(T)\cup V(T^{\prime})\setminus\{v,v^{\prime}\}\}|\,.$$

Let $w \in V(G_0)$. Since $G_0$ is a $k$-$(v,v^{\prime})$-UD graph, we have $\ecc_G(w)=\ecc_{G_0}(w)+r$ and hence
\begin{align*}
\tr_G(w) & = \sum\limits_{z\in V(G_0)}d_G(w,z)+\sum\limits_{z\in V(T)\setminus\{v\}}d_G(w,z)+\sum\limits_{z\in V(T^{\prime})\setminus\{v^{\prime}\}}d_G(w,z) \\
& = \tr_{G_0}(w)+\sum\limits_{z\in V(T)\setminus\{v\}}\Big[d_{G_0}(w,v)+d_T(v,z)\Big] \\
&\quad + \sum\limits_{z\in V(T^{\prime})\setminus\{v^{\prime}\}}\Big[d_{G_0}(w,v^{\prime})+d_{T^{\prime}}(v^{\prime},z)\Big] \\
&= \tr_{G_0}(w)+(n_1-1)\Big[d_{G_0}(w,v)+d_{G_0}(w,v^{\prime})\Big]+\tr_{T}(v)+\tr_{T^{\prime}}(v^{\prime})\\
&=\tr_{G_0}(w)+(n_1-1)k+2\tr_{T}(v)
\end{align*} for every vertex $w\in V(G_0)$. Thus $|\{\tr_G(w):w\in V(G_0)\}|=|\{\tr_{G_0}(w):w\in V(G_0)\}|$ and $|\{\ecc_G(w):w\in V(G_0)\}|=|\{\ecc_{G_0}(w):w\in V(G_0)\}|$. Notice that as $k = \diam(G_0)$, $\ecc_G(y) = r + k + i > \ecc_{G_0}(w) + r = \ecc_G(w)$ for every $y \in V(T) \setminus \{v\}$ and every $w \in V(G_0)$, thus $C_{\ecc}(G)=C_{\ecc}(G_0)+r$. Since $C_{\ecc}(G_0)>C_W(G_0)$, we conclude that $C_{\ecc}(G)=C_{\ecc}(G_0)+r>C_W(G_0)+r\geq C_W(G)$.
\qed

Denote by $Y_k$ the graph of order $2k+4$  consisting of $P_{2k+3}$ and an additional vertex which is adjacent to the two neighbors of the  central vertex of $P_{2k+3}$. Note that $Y_1=Z_1$  with $C_{\ecc}(Z_1)>C_W(Z_1)$ and $Y_k$ with $k>1$ is obtained by attaching a pendant path of $k-1$ vertices to each pedant vertex of $Z_1$. By similar reasoning as that in the proof of Theorem~\ref{New-1}, we have $C_{\ecc}(Y_k)>C_W(Y_k)$.

Now we have the following  natural question. Is it true that if $C_{\ecc}(G)>C_W(G)$ and $C_{\ecc}(H)>C_W(H)$, then $C_{\ecc}(G\cp H)>C_W(G\cp H)$ holds? The answer is negative, as can be seen from the following result: $$C_{\ecc}(Y_k\cp Z_k)-C_{W}(Y_k\cp Z_k)=\left\{
\begin{array}{rl}
1, & k=2\,; \\
0,& k\in \{3,4\}\,; \\
-2,& k=5\,.
\end{array}
\right.$$

The graphs $G$ with $C_{\ecc}(G) > C_W(G)$ found in this paper have diameter at least $4$. On the other hand, it follows from the proof of Proposition~\ref{prp:for-almosst-all} that there are no such graphs of diameter $2$. Hence, we pose:

\begin{problem}
	Does there exist a graph $G$ with $\diam(G) = 3$ and $C_{\ecc}(G) > C_W(G)$?
\end{problem}

 We have checked by computer that there is no such graph of order at most $10$.
Since $1\leq C_{\ecc}(G)\leq 2$ for any graph $G$ with $\diam(G)=3$, the key point for solving the above problem is to determine the existence of transmission regular but non-self-centered graphs with diameter $3$.

\section*{Acknowledgements}
 The authors thank the anonymous referees for their strict criticisms and helpful suggestions which improve the presentation of this paper.
Kexiang Xu is supported by supported by NNSF of China (grant No.\ 11671202, and the China-Slovene bilateral grant 12-9). Sandi Klav\v{z}ar acknowledges the financial support from the Slovenian Research Agency (research core funding P1-0297, projects J1-9109, J1-1693, N1-0095, and the bilateral grant BI-CN-18-20-008).


\end{document}